\documentclass[12pt]{article}
\usepackage{mathrsfs}
\usepackage{amsfonts}
\usepackage{amsmath}
\usepackage{psfrag}
\usepackage{amsfonts}
\makeatletter
\newcommand{\rmnum}[1]{\romannumeral #1}
\newcommand{\Rmnum}[1]{\expandafter\@slowromancap\romannumeral #1@}
\makeatother
\title{Rigidity of Gradient Shrinking Ricci Solitons\footnote{This work is partially supported by Natural Science Foundation of China (No. 11601495).}}
\author{\small Fei Yang\footnote{ School of Mathematics and Physics, China University of Geosciences, Wuhan, People's Republic of China. E-mail address: yangfei810712@163.com.}\and \small  Liangdi Zhang\footnote{ School of Mathematics and Physics, China University of Geosciences, Wuhan, People's Republic of China. E-mail address: ldzhang91@163.com.}}
\date{}
\begin{document}
\maketitle{}
\begin{abstract}
We prove that a gradient shrinking Ricci soliton with fourth order divergence-free Riemannian tensor is rigid. For the $4$-dimensional case, we show that any gradient shrinking Ricci soliton with fourth order divergence-free Riemannian tensor is either Einstein, or a finite quotient of the Gaussian shrinking soliton
$\mathbb{R}^4$, $\mathbb{R}^2\times\mathbb{S}^2$ or the round cylinder $\mathbb{R}\times\mathbb{S}^3$. Under the condition of fourth order divergence-free Weyl tensor, we have the same results.
\newline
\end{abstract}
\section{Introduction}
\renewcommand{\theequation}{\thesection.\arabic{equation}}
\indent\par
A complete Riemannian manifold $(M^n,g,f)$ is called a gradient Ricci soliton if there exists a smooth function $f$ on $M^n$ such that the Ricci tensor $Ric$ of the metric $g$ satisfies the equation
\begin{equation}
Ric+\nabla^2f=\lambda g
\end{equation}
for some constant $\lambda$. For $\lambda>0$ the Ricci soliton is shrinking, for $\lambda=0$ it is steady and for $\lambda<0$ expanding.

The classification of gradient Ricci solitons has been a subject of interest for many people in recent years. For four-dimensional gradient Ricci solitons, A. Naber [12] showed that a four-dimensional non-compact shrinking Ricci soliton with bounded nonnegative Riemannian curvature is a finite quotient of $\mathbb{R}^4$, $\mathbb{R}^2\times\mathbb{S}^2$ or $\mathbb{R}\times\mathbb{S}^3$. X. Chen and Y. Wang [6] classified four-dimensional anti-self dual gradient steady and shrinking Ricci solitons. More generally, J.Y. Wu, P. Wu and W. Wylie [17] proved that a four-dimensional gradient shrinking Ricci soliton with half harmonic Weyl tensor (i.e. $divW^\pm=0$) is either Einstein or a finite quotient of $\mathbb{R}^4$, $\mathbb{R}^2\times\mathbb{S}^2$ or $\mathbb{R}\times\mathbb{S}^3$.

For $n$-dimensional gradient Ricci solitons, M. Eminenti, G. La Nave and C. Mantegazza [7] proved that an $n$-dimensional compact shrinking Ricci soliton with vanishing Weyl tensor is a finite quotient of $\mathbb{S}^n$. More generally, P. Peterson and W. Wylie [14] showed that a gradient shrinking Ricci soliton with vanishing Weyl tensor is a finite quotient of $\mathbb{R}^n$, $\mathbb{S}^{n-1}\times\mathbb{R}$, or $\mathbb{S}^n$ by assuming $\int_M|Ric|^2e^{-f}<\infty$. The integral assumption was proven to be true for gradient shrinking Ricci solitions (see Theorem 1.1 of [11]). Without additional assumptions, Z. H. Zhang [18] obtained the same classification of gradient shrinking Ricci solitons with vanishing Weyl tensor.

H. D. Cao and Q. Chen [1] introduced the covariant 3-tensor $D$, i.e.
\begin{eqnarray*}
D_{ijk}&=&\frac{1}{n-2}(R_{jk}\nabla_if-R_{ik}\nabla_jf)+\frac{1}{2(n-1)(n-2)}(g_{jk}\nabla_iR-g_{ik}\nabla_jR)\\
&&-\frac{R}{(n-1)(n-2)}(g_{jk}\nabla_if-g_{ik}\nabla_jf).
\end{eqnarray*}
to study the classification of locally conformally flat gradient steady
solitons. The vanishing of $D$ is a crucial ingredient in their classification results. They [2] proved that a compact gradient shrinking Ricci solitons with $D=0$ is Einstein. Moreover, they showed that a four-dimensional complete non-compact Bach-flat gradient shrinking Ricci soliton is a finite quotient of $\mathbb{R}^4$ or $\mathbb{R}\times\mathbb{S}^3$. More generally, they proved that a $n$-dimensional $(n\geq5)$ complete non-compact Bach-flat gradient shrinking Ricci soliton with is a finite quotient of $\mathbb{R}^n$ or $\mathbb{R}\times N^{n-1}$, where $N$ is an $(n-1)$-dimensional Einstein manifold.

 H. D. Cao and Q. Chen [2] proved that a Bach-flat gradient shrinking Ricci soliton has vanishing $D$, where the Bach tensor is given by
$$B_{ij}=\frac{1}{n-3}\nabla_k\nabla_lW_{ikjl}+\frac{1}{n-2}R_{kl}W_{ikjl}.$$
 Moreover, they showed that the $3$-tensor $D$ is closely related to Cotton tensor, i.e.
 $$C_{ijk}=\nabla_iR_{jk}-\nabla_jR_{ik}-\frac{1}{2(n-1)}(g_{jk}\nabla_iR-g_{ik}\nabla_jR),$$
and Weyl tensor, i.e.
\begin{eqnarray*}
W_{ijkl}&=&R_{ijkl}-\frac{1}{n-2}(g_{ik}R_{jl}-g_{il}R_{jk}-g_{jk}R_{il}+g_{jl}R_{ik})\\
&&+\frac{R}{(n-1)(n-2)}(g_{ik}g_{jl}-g_{il}g_{jk})
\end{eqnarray*}
by
\[D_{ijk}=C_{ijk}+W_{ijkl}\nabla_l f.\]

M. Fern\'{a}ndez-L\'{o}pez and E. Garcia-R\'{i}o [8] proved that the compact Ricci soliton is rigid if and only if it has harmonic Weyl tensor. For the complete non-compact case, O. Munteanu and N. Sesum [11] showed that a gradient shrinking Ricci soliton with harmonic Weyl tensor is rigid.

In 2016, G. Catino, P. Mastrolia and D. D. Monticelli [4] proved that the gradient shrinking Ricci soliton is rigid if $div^4W=0$. In their paper, $div^4$ is defined by $div^4W=\nabla_k\nabla_j\nabla_l\nabla_iW_{ikjl}$. They showed that $div^4W=0$ if and only if $div^3C=0$, where $div^3C=\nabla_i\nabla_j\nabla_kC_{ijk}$. Then, they proved that $div^3C=0$ implies $C=0$. The rigidity result follows.

S. Tachibana [15] proved that a compact orientable Riemannian manifold with $Rm>0$ and $divRm=0$ is a space of constant curvature. P. Peterson and W. Wylie [13] proved that a compact shrinking gradient Ricci soliton is Einstein if $\int_M Ric(\nabla f,\nabla f)\leq0$. They also showed that a gradient Ricci soliton is rigid if and only if it has constant scalar curvature and is radially flat.

In order to  state our results precisely, we introduce the following definitions for the Riemannian curvature:
\[(divRm)_{ijk}:=\nabla_lR_{ijkl},\]
\[(div^2Rm)_{ik}:=\nabla_j\nabla_lR_{ijkl},\]
\[(div^3Rm)_i:=\nabla_k\nabla_j\nabla_lR_{ijkl},\]
\[div^4Rm:=\nabla_i\nabla_k\nabla_j\nabla_lR_{ijkl}.\]

For the Weyl curvature tensor, we define:
\[(divW)_{ijk}:=\nabla_lW_{ijkl},\]
\[(div^2W)_{ik}:=\nabla_j\nabla_lW_{ijkl},\]
\[(div^3W)_i:=\nabla_k\nabla_j\nabla_lW_{ijkl},\]
\[div^4W:=\nabla_i\nabla_k\nabla_j\nabla_lW_{ijkl}.\]

Our main results are the following theorems for gradient shrinking Ricci solitons:
\\\\\textbf{Theorem 1.1} Let $(M^n,g)$ be a gradient shrinking Ricci soliton with (1.1). If $div^4Rm=0$, then $(M^n,g)$ is rigid.
\\\\\textbf{Theorem 1.2} Let $(M^n,g)$ be a gradient shrinking Ricci soliton with (1.1). If $div^3Rm(\nabla f)=0$, then $(M^n,g)$ is rigid.

G. Catino, P. Mastrolia and D. D. Monticelli [4] proved that a gradient shrinking Ricci soliton with $div^4W=0$ is rigid. We will give a different proof in Section 8 Appendix. Moreover, we have the following result:
\\\\\textbf{Theorem 1.3} Let $(M^n,g)$ be a gradient shrinking Ricci soliton with (1.1). If $div^3W(\nabla f)=0$, then $(M^n,g)$ is rigid.

For the $4$-dimensional case, we have the following classification theorems:
\\\\\textbf{Theorem 1.4} Let $(M^4,g)$ be a $4$-dimensional gradient shrinking Ricci soliton with (1.1). If $div^4Rm=0$, then $(M^4,g)$ is either

(\rmnum{1}) Einstein, or

(\rmnum{2}) a finite quotient of the Gaussian shrinking soliton
$\mathbb{R}^4$, $\mathbb{R}^2\times\mathbb{S}^2$ or the round cylinder $\mathbb{R}\times\mathbb{S}^3$.
\\\\\textbf{Theorem 1.5} Let $(M^4,g)$ be a $4$-dimensional gradient shrinking Ricci soliton with (1.1). If $div^3Rm(\nabla f)=0$,  then $(M^4,g)$ is either

(\rmnum{1}) Einstein, or

(\rmnum{2}) a finite quotient of the Gaussian shrinking soliton
$\mathbb{R}^4$, $\mathbb{R}^2\times\mathbb{S}^2$ or the round cylinder $\mathbb{R}\times\mathbb{S}^3$.
\\\\\textbf{Theorem 1.6} Let $(M^4,g)$ be a $4$-dimensional gradient shrinking Ricci soliton with (1.1). If $div^3W(\nabla f)=0$,  then $(M^4,g)$ is either

(\rmnum{1}) Einstein, or

(\rmnum{2}) a finite quotient of the Gaussian shrinking soliton
$\mathbb{R}^4$, $\mathbb{R}^2\times\mathbb{S}^2$ or the round cylinder $\mathbb{R}\times\mathbb{S}^3$.
\\\\\textbf{Remark 1.1} As it will be clear from the proof, the scalar assumptions on the vanishing
of $div^4Rm$, $div^3Rm(\nabla f)$, and $div^3W(\nabla f)$ in all the above theorems can be trivially relaxed to a (suitable) inequality. To be precise, Theorem 1.1 and Theorem 1.4 hold just assuming $div^4Rm\geq0$. Under the condition of $div^3Rm(\nabla f)\geq0$, Theorem 1.2 and Theorem 1.5 still hold. Moreover, Theorem 1.3 and Theorem 1.6 hold for $div^3W(\nabla f)\geq0$.

The rest of this paper is organized as follows. In Section 2, we recall some background material and prove some formulas which will be needed in the proof of the main theorems. In Section 3, we will prove that the compact gradient Ricci soliton with fourth divergence-free Riemannian tensor is Einstein. The proof makes use of a rigid theorem obtained by P. Peterson and W. Wylie [13]. In Section 4, we will deal with the complete noncompact case of Theorem 1.1. In Section 5, we give a direct proof of Theorem 1.2. We first prove divergence formulas of the Weyl tensor in Section 6, then we will prove Theorem 1.3. Finally, in Section 7 we will finish the proof of Theorems 1.4 to 1.6.
\section{Preliminaries}
\renewcommand{\theequation}{\thesection.\arabic{equation}}
\indent\par
First of all, we present some basic facts of gradient shrinking Ricci solitons.
\\\\\textbf{Proposition 2.1} (([7,10,11,14])) Let $(M^n,g)$ be a gradient shrinking Ricci soliton with (1.1), we have the following identities.
\begin{equation}
\nabla_lR_{ijkl}=\nabla_jR_{ik}-\nabla_iR_{jk},
\end{equation}
\begin{equation}
\nabla R=2divRic,
\end{equation}
\begin{equation}
R_{ijkl}\nabla_lf=\nabla_lR_{ijkl},
\end{equation}
\begin{equation}
\nabla_l(R_{ijkl}e^{-f})=0,
\end{equation}
\begin{equation}
R_{jl}\nabla_lf=\nabla_lR_{jl},
\end{equation}
\begin{equation}
\nabla_l(R_{jl}e^{-f})=0,
\end{equation}
\begin{equation}
\nabla R=2Ric(\nabla f,\cdot),
\end{equation}
\begin{equation}
\Delta_fR_{ik}=2\lambda R_{ik}-2R_{ijkl}R_{jl},
\end{equation}
\begin{equation}
\Delta_fR=2\lambda R-2|Ric|^2,
\end{equation}
where $\Delta_f:=\Delta-\nabla_{\nabla f}$,
\begin{equation}
\Delta_f|Ric|^2=4\lambda|Ric|^2-4Rm(Ric,Ric)+2|\nabla Ric|^2,
\end{equation}
where $Rm(Ric,Ric)=R_{ijkl}R_{ik}R_{jl}$, and
\begin{equation}
R+|\nabla f|^2-2\lambda f=Const.
\end{equation}

Next we prove the following formulas for gradient shrinking Ricci soliton with (1.1).
\\\\\textbf{Proposition 2.2} Let $(M^n,g)$ be a gradient shrinking Ricci soliton with (1.1), we have the following identities.
\begin{equation}
(div^2Rm)_{ik}=2\lambda R_{ik}+\nabla_lR_{ik}\nabla_lf-\frac{1}{2}\nabla_i\nabla_kR-R_{ik}^2-R_{ijkl}R_{jl},
\end{equation}
\begin{equation}
(div^3Rm)_{i}=-R_{ijkl}\nabla_kR_{jl},
\end{equation}
and
\begin{equation}
div^4Rm=\nabla_lR_{jk}\nabla_kR_{jl}-|\nabla Ric|^2-R_{ijkl}\nabla_i\nabla_kR_{jl}.
\end{equation}
\\\\\textbf{Proof.} By direct computation,
\begin{eqnarray*}
(div^2Rm)_{ik}&=&\nabla_j\nabla_lR_{ijkl}\notag\\
&=&\Delta R_{ik}-\nabla_j\nabla_iR_{jk}\notag\\
&=&\Delta_fR_{ik}+\nabla_lR_{ik}\nabla_lf-\nabla_i\nabla_jR_{jk}+R_{ijkl}R_{jl}-R_{ik}^2\notag\\
&=&2\lambda R_{ik}-2R_{ijkl}R_{jl}+\nabla_lR_{ik}\nabla_lf-\frac{1}{2}\nabla_i\nabla_kR+R_{ijkl}R_{jl}-R_{ik}^2\notag\\
&=&2\lambda R_{ik}-R_{ijkl}R_{jl}+\nabla_lR_{ik}\nabla_lf-\frac{1}{2}\nabla_i\nabla_kR-R_{ik}^2,
\end{eqnarray*}
where we used (2.2) in the second equality. Moreover, we used (2.3) and (2.9) in the fourth equality.

Using (2.13), we have
\begin{eqnarray*}
(div^3Rm)_{i}&=&\nabla_k\nabla_j\nabla_lR_{ijkl}\notag\\
&=&\nabla_k(2\lambda R_{ik}-R_{ijkl}R_{jl}+\nabla_lR_{ik}\nabla_lf-\frac{1}{2}\nabla_i\nabla_kR-R_{ik}^2)\notag\\
&=&\lambda\nabla_iR-\nabla_kR_{ijkl}R_{jl}-R_{ijkl}\nabla_kR_{jl}+\nabla_lR_{ik}\nabla_k\nabla_lf+\nabla_k\nabla_lR_{ik}\nabla_lf\notag\\
&&-\frac{1}{2}\nabla_k\nabla_i\nabla_kR-R_{ij}\nabla_kR_{kj}-R_{kj}\nabla_kR_{ij}\notag\\
&=&\lambda\nabla_iR+(\nabla_jR_{il}-\nabla_iR_{jl})R_{jl}-R_{ijkl}\nabla_kR_{jl}+\nabla_lR_{ik}(\lambda g_{kl}-R_{kl})\notag\\
&&+(\nabla_l\nabla_kR_{ik}+R_{lj}R_{ij}+R_{klij}R_{jk})\nabla_lf-\frac{1}{2}\nabla_i\Delta_fR-\frac{1}{2}\nabla_i(\nabla_kR\nabla_kf)\notag\\
&&-\frac{1}{2}R_{ij}\nabla_jR-\frac{1}{2}R_{ij}\nabla_jR-R_{kj}\nabla_kR_{ij}\notag\\
&=&\lambda\nabla_iR+R_{jl}\nabla_jR_{il}-\frac{1}{2}\nabla_i|Ric|^2-R_{ijkl}\nabla_kR_{jl}+\frac{\lambda}{2}\nabla_iR\notag\\
&&-R_{kl}\nabla_lR_{ik}+\frac{1}{2}\nabla_l\nabla_iR\nabla_lf+\frac{1}{2}R_{ij}\nabla_jR+R_{jk}\nabla_lR_{ijkl}\notag\\
&&-\lambda\nabla_iR+\nabla_i|Ric|^2-\frac{1}{2}\nabla_i\nabla_lR\nabla_lf-\frac{1}{2}\nabla_lR\nabla_i\nabla_lf\notag\\
&&-R_{ij}\nabla_jR-R_{kj}\nabla_kR_{ij}\notag\\
&=&\frac{1}{2}\nabla_i|Ric|^2-R_{ijkl}\nabla_kR_{jl}+\frac{\lambda}{2}\nabla_iR-\frac{1}{2}R_{ik}\nabla_kR+R_{jk}\nabla_jR_{ik}\notag\\
&&-\frac{1}{2}\nabla_i|Ric|^2-\frac{\lambda}{2}\nabla_iR+\frac{1}{2}R_{il}\nabla_lR-R_{kj}\nabla_kR_{ij}\notag\\
&=&-R_{ijkl}\nabla_kR_{jl},
\end{eqnarray*}
where we used (2.3) in the third equality, used (2.2) and (1.1) in the fourth equality. Moreover, we used (2.4), (2.8) and (2.10) in the fifth equality. In the sixth equality, we used  (1.1) and (2.2).

It follows from (2.14) that
\begin{eqnarray*}
div^4Rm&=&\nabla_i\nabla_k\nabla_j\nabla_lR_{ijkl}\notag\\
&=&-\nabla_iR_{ijkl}\nabla_kR_{jl}-R_{ijkl}\nabla_i\nabla_kR_{jl}\notag\\
&=&(\nabla_lR_{jk}-\nabla_kR_{jl})\nabla_kR_{jl}-R_{ijkl}\nabla_i\nabla_kR_{jl}\notag\\
&=&\nabla_lR_{jk}\nabla_kR_{jl}-|\nabla Ric|^2-R_{ijkl}\nabla_i\nabla_kR_{jl},
\end{eqnarray*}
where we used (2.2) in the third equality.$\hfill\Box$
\\\\\textbf{Remark 2.1} It follows from (2.13) that $div^2Rm$ is a symmetric 2-tensor. Therefore, we have the following identities.
\[(div^2Rm)_{ik}=\nabla_j\nabla_lR_{ijkl}=\nabla_l\nabla_jR_{ijkl},\]
\[(div^3Rm)_{i}=\nabla_k\nabla_j\nabla_lR_{ijkl}=\nabla_k\nabla_j\nabla_lR_{kjil}=\nabla_k\nabla_l\nabla_jR_{ijkl}=\nabla_k\nabla_l\nabla_jR_{kjil},\]
and
\[div^4Rm=\nabla_i\nabla_k\nabla_j\nabla_lR_{ijkl}=\nabla_i\nabla_k\nabla_l\nabla_jR_{ijkl}=\nabla_k\nabla_i\nabla_j\nabla_lR_{ijkl}=\nabla_k\nabla_i\nabla_l\nabla_jR_{ijkl}.\]

Finally, we list following results that will be needed in the proof of the main theorems.
\\\\\textbf{Lemma 2.1}  Let $(M^n,g)$ be a complete gradient shrinking soliton with (1.1). Then it has
nonnegative scalar curvature $R\geq0$.
\\\\\textbf{Remark 2.2} Lemma 2.1 is a special case of a more general result of B. L. Chen [5] which
states that $R\geq0$ for any ancient solution to the Ricci flow.
\\\\\textbf{Lemma 2.2} (H. D. Cao and D. Zhou [3])  Let $(M^n,g)$ be a complete gradient shrinking soliton with (1.1). Then,

(\rmnum{1}) the potential function f satisfies the estimates
\begin{equation}
\frac{1}{4}(r(x)-c_1)^2\leq f(x)\leq \frac{1}{4}(r(x)+c_2)^2,
\end{equation}
where $r(x)=d(x_0,x)$ is the distance function from some fixed point $x_0\in M$, $c_1$
and $c_2$ are positive constants depending only on $n$ and the geometry of $g$ on the
unit ball $B(x_0,1)$;

(\rmnum{2}) there exists some constant $C>0$ such that
\begin{equation}
Vol(B(x_0,s))\leq Cs^n
\end{equation}
for $s>0$ sufficiently large.
\\\\\textbf{Lemma 2.3} (P. Petersen and W. Wylie [13]) A shrinking compact gradient soliton is rigid with trivial $f$ if
\begin{equation}
\int_MRic(\nabla f,\nabla f)\leq0.
\end{equation}
\\\\\textbf{Lemma 2.4} ( P. Petersen and W. Wylie [13]) A gradient soliton is rigid if and only if it has
constant scalar curvature and is radially flat, that is, $sec(E,\nabla f)=0$.
\\\\\textbf{Remark 2.3} The condition of $divRm=0$ is stronger than $sec(E,\nabla f)=0$.
\\\\\textbf{Lemma 2.5} (O. Munteanu and N. Sesum [11]) For any complete gradient shrinking Ricci soliton with (1), we have
\begin{equation}
\int_M|Ric|^2e^{-\alpha f}<+\infty
\end{equation}
for any $\alpha>0$.
\\\\\textbf{Lemma 2.6} (O. Munteanu and N. Sesum [11]) Let $(M, g)$ be a gradient shrinking Ricci soliton. If for some $\beta<1$ we
have $\int_M|Rm|^2e^{-\beta f}<+\infty$, then the following identity holds.
\begin{equation}
\int_M|div Rm|^2e^{-f}=\int_M|\nabla Ric|^2e^{-f}<+\infty.
\end{equation}
\section{The Compact Case of Theorem 1.1}
\renewcommand{\theequation}{\thesection.\arabic{equation}}
\indent\par
In this section, we prove the the compact case of Theorem 1.1:
\\\\\textbf{Theorem 3.1} Let $(M^n,g)$ be a compact gradient shrinking Ricci soliton with (1.1). If $div^4Rm=0$, then $(M^n,g)$ is Einstein.

The first step in proving Theorem 3.1 is to obtain the following integral equation.
\\\\\textbf{Lemma 3.1} Let $(M^n,g)$ be a compact gradient shrinking Ricci soliton with (1.1), then
\begin{equation}
\int_M\nabla_lR_{jk}\nabla_kR_{jl}e^{-f}=\frac{1}{2}\int_M|\nabla Ric|^2e^{-f}.
\end{equation}
\\\\\textbf{Proof.} Calculating directly, we have
\begin{eqnarray}
&&\int_M\nabla_lR_{jk}\nabla_kR_{jl}e^{-f}\notag\\
&=&-\int_MR_{jk}\nabla_l\nabla_kR_{jl}e^{-f}+\int_MR_{jk}\nabla_kR_{jl}\nabla_lfe^{-f}\notag\\
&=&-\int_MR_{jk}(\nabla_k\nabla_lR_{jl}+R_{jp}R_{pk}+R_{lkji}R_{il})e^{-f}\notag\\
&&+\int_MR_{jk}\nabla_k(R_{jl}\nabla_lf)e^{-f}-\int_MR_{jk}R_{jl}\nabla_k\nabla_lfe^{-f}\notag\\
&=&-\int_MR_{jk}(R_{jp}R_{pk}+R_{lkji}R_{il})e^{-f}-\int_MR_{jk}R_{jl}(\lambda g_{kl}-R_{kl})e^{-f}\notag\\
&=&-\int_MtrRic^3e^{-f}+\int_MRm(Ric,Ric)e^{-f}-\lambda\int_M|Ric|^2e^{-f}+\int_MtrRic^3e^{-f}\notag\\
&=&\int_MRm(Ric,Ric)e^{-f}-\lambda\int_M|Ric|^2e^{-f},
\end{eqnarray}
where we used (2.6) and (1.1) in the third equality.

Applying (2.11) to (3.22), we obtain
\begin{eqnarray*}
&&\int_M\nabla_lR_{jk}\nabla_kR_{jl}e^{-f}\notag\\
&=&-\frac{1}{4}\int_M\Delta_f|Ric|^2e^{-f}+\frac{1}{2}\int_M|\nabla Ric|^2e^{-f}\notag\\
&=&-\frac{1}{4}\int_M(\Delta|Ric|^2-\nabla_{\nabla f}|Ric|^2)e^{-f}+\frac{1}{2}\int_M|\nabla Ric|^2e^{-f}\notag\\
&=&-\frac{1}{4}\int_M\nabla_{\nabla f}|Ric|^2e^{-f}+\frac{1}{4}\int_M\nabla_{\nabla f}|Ric|^2e^{-f}+\frac{1}{2}\int_M|\nabla Ric|^2e^{-f}\notag\\
&=&\frac{1}{2}\int_M|\nabla Ric|^2e^{-f}.
\end{eqnarray*}
$\hfill\Box$

Now we are ready to prove Theorem 3.1.
\\\\\textbf{Proof of Theorem 3.1:}

Integrating (2.15), we obtain
\begin{eqnarray}
&&\int_Mdiv^4Rme^{-f}\notag\\
&=&\int_M\nabla_lR_{jk}\nabla_kR_{jl}e^{-f}-\int_M|\nabla Ric|^2e^{-f}-\int_MR_{ijkl}\nabla_i\nabla_kR_{jl}e^{-f}\notag\\
&=&\frac{1}{2}\int_M|\nabla Ric|^2e^{-f}-\int_M|\nabla Ric|^2e^{-f}\notag\\
&=&-\frac{1}{2}\int_M|\nabla Ric|^2e^{-f},
\end{eqnarray}
where we used Lemma 3.1 and (2.5) in the second equality.

Since $div^4Rm=0$, it follows from (3.23) that $\int_M|\nabla Ric|^2e^{-f}=0$, i.e. $|\nabla Ric|=0$ $a. e$. . Note that any gradient shrinking Ricci soliton is analytic in harmonic coordinates, we have $|\nabla Ric|=0$ on $M$.

By direct computation, we have
\[0\leq|\nabla Ric-\frac{\nabla R}{n}g|^2=|\nabla Ric|^2-\frac{|\nabla R|^2}{n}=-\frac{|\nabla R|^2}{n}\leq0.\]

Therefore, $R$ is a constant on $M$. It follows from (2.8) that $Ric(\nabla f,\nabla f)=\frac{1}{2}\langle\nabla R,\nabla f\rangle=0$. By Lemma 2.3, $(M^n,g)$ is rigid. The compactness of $(M^n,g)$ implies that $(M^n,g)$ is Einstein.$\hfill\Box$
\section{The Complete Non-compact Case of Theorem 1.1}
\renewcommand{\theequation}{\thesection.\arabic{equation}}
\indent\par
In this section, we prove the complete non-compact case of Theorem 1.1:
\\\\\textbf{Theorem 4.1} Let $(M^n,g)$ be a complete non-compact gradient shrinking Ricci soliton with (1.1). If $div^4Rm=0$, then $(M^n,g)$ is rigid.

The first step in proving Theorem 4.1 is to obtain the following integral inequality.
\\\\\textbf{Lemma 4.1} Let $(M^n,g)$ be a complete non-compact gradient shrinking Ricci soliton with (1.1). For every $C^2$ function $\phi:\mathbb{R_+}\rightarrow\mathbb{R}$ with $\phi(f)$ having compact support in $M$ and some constant $c>0$, we have
\begin{equation}
\int_M\nabla_lR_{jk}\nabla_kR_{jl}\phi^2(f)e^{-f}\leq c\int_M|Ric|^2|\nabla f|^2(\phi')^2e^{-f}+\frac{3}{4}\int_M|\nabla Ric|^2\phi^2(f)e^{-f}.
\end{equation}
\\\\\textbf{Proof.} By direct computation, we have
\begin{eqnarray}
&&\int_M\nabla_lR_{jk}\nabla_kR_{jl}\phi^2(f)e^{-f}\notag\\
&=&-\int_MR_{jk}\nabla_l\nabla_kR_{jl}\phi^2(f)e^{-f}-\int_MR_{jk}\nabla_kR_{jl}\nabla_l\phi^2(f)e^{-f}+\int_MR_{jk}\nabla_kR_{jl}\nabla_lf\phi^2(f)e^{-f}\notag\\
&=&-\int_MR_{jk}(\nabla_k\nabla_lR_{jl}+R_{jp}R_{pk}+R_{lkji}R_{il})\phi^2(f)e^{-f}-2\int_MR_{jk}\nabla_kR_{jl}\nabla_lf\phi\phi'e^{-f}\notag\\
&&+\int_MR_{jk}\nabla_k(R_{jl}\nabla_lf)\phi^2(f)e^{-f}-\int_MR_{jk}R_{jl}\nabla_k\nabla_lf\phi^2(f)e^{-f}\notag\\
&=&-\int_MR_{jk}(R_{jp}R_{pk}+R_{lkji}R_{il})\phi^2(f)e^{-f}-2\int_MR_{jk}\nabla_kR_{jl}\nabla_lf\phi\phi'e^{-f}\notag\\
&&-\int_MR_{jk}R_{jl}(\lambda g_{kl}-R_{kl})\phi^2(f)e^{-f}\notag\\
&=&-\int_MtrRic^3\phi^2(f)e^{-f}+\int_MRm(Ric,Ric)\phi^2(f)e^{-f}-2\int_MR_{jk}\nabla_kR_{jl}\nabla_lf\phi\phi'e^{-f}\notag\\
&&-\lambda\int_M|Ric|^2\phi^2(f)e^{-f}+\int_MtrRic^3\phi^2(f)e^{-f}\notag\\
&=&\int_MRm(Ric,Ric)\phi^2(f)e^{-f}-2\int_MR_{jk}\nabla_kR_{jl}\nabla_lf\phi\phi'e^{-f}-\lambda\int_M|Ric|^2\phi^2(f)e^{-f},\notag\\
\end{eqnarray}
where we used (2.7) and (1.1) in the third equality.

Applying (2.11) to (4.25), we obtain
\begin{eqnarray*}
&&\int_M\nabla_lR_{jk}\nabla_kR_{jl}\phi^2(f)e^{-f}\notag\\
&=&-2\int_MR_{jk}\nabla_kR_{jl}\nabla_lf\phi\phi'e^{-f}-\frac{1}{4}\int_M\Delta_f|Ric|^2\phi^2(f)e^{-f}+\frac{1}{2}\int_M|\nabla Ric|^2\phi^2(f)e^{-f}\notag\\
&=&-2\int_MR_{jk}\nabla_kR_{jl}\nabla_lf\phi\phi'e^{-f}-\frac{1}{4}\int_M\Delta|Ric|^2\phi^2(f)e^{-f}\notag\\
&&+\frac{1}{4}\int_M\nabla_{\nabla f}|Ric|^2\phi^2(f)e^{-f}+\frac{1}{2}\int_M|\nabla Ric|^2\phi^2(f)e^{-f}\notag\\
&=&-2\int_MR_{jk}\nabla_kR_{jl}\nabla_lf\phi\phi'e^{-f}+\frac{1}{4}\int_M\langle\nabla|Ric|^2,\nabla\phi^2(f)\rangle e^{-f}\notag\\
&&-\frac{1}{4}\int_M\nabla_{\nabla f}|Ric|^2\phi^2(f)e^{-f}+\frac{1}{4}\int_M\nabla_{\nabla f}|Ric|^2\phi^2(f)e^{-f}+\frac{1}{2}\int_M|\nabla Ric|^2\phi^2(f)e^{-f}\notag\\
&=&-2\int_MR_{jk}\nabla_kR_{jl}\nabla_lf\phi\phi'e^{-f}+\int_MR_{ik}\nabla_lR_{ik}\nabla_lf\phi\phi'e^{-f}\notag\\
&&+\frac{1}{2}\int_M|\nabla Ric|^2\phi^2(f)e^{-f}\notag\\
&\leq&c\int_M|Ric||\nabla f||\nabla Ric||\phi||\phi'|e^{-f}+\frac{1}{2}\int_M|\nabla Ric|^2\phi^2(f)e^{-f}\notag\\
&\leq&c\int_M|Ric|^2|\nabla f|^2(\phi')^2e^{-f}+\frac{3}{4}\int_M|\nabla Ric|^2\phi^2(f)e^{-f}
\end{eqnarray*}
for some constant $c>0$.$\hfill\Box$
\\\\\textbf{Lemma 4.2} Let $(M^n,g)$ be a complete non-compact gradient shrinking Ricci soliton with (1.1). For every $C^2$ function $\varphi:\mathbb{R_+}\rightarrow\mathbb{R}$ with $\varphi(f)$ having compact support in $M$, we have
\begin{equation}
-\int_MR_{ijkl}\nabla_i\nabla_kR_{jl}\varphi(f)e^{-f}=\int_M(|\nabla Ric|^2-\nabla_lR_{kj}\nabla_kR_{jl})\varphi'e^{-f}.
\end{equation}
\\\\\textbf{Proof.} By direct computation, we have
\begin{eqnarray*}
-\int_MR_{ijkl}\nabla_i\nabla_kR_{jl}\varphi(f)e^{-f}&=&\int_MR_{ijkl}\nabla_kR_{jl}\varphi'\nabla_ife^{-f}\notag\\
&=&\int_M\nabla_iR_{ijkl}\nabla_kR_{jl}\varphi'e^{-f}\notag\\
&=&\int_M(\nabla_kR_{jl}-\nabla_lR_{kj})\nabla_kR_{jl}\varphi'e^{-f}\notag\\
&=&\int_M(|\nabla Ric|^2-\nabla_lR_{kj}\nabla_kR_{jl})\varphi'e^{-f},
\end{eqnarray*}
where we used (2.5), (2.4) and (2.2) in the first, second and third equality, respectively.$\hfill\Box$

Now we are ready to prove Theorem 4.1.
\\\\\textbf{Proof of Theorem 4.1:}

Let $\phi:\mathbb{R_+}\rightarrow\mathbb{R}$ be a $C^2$ function with $\phi=1$ on $(0,s]$, $\phi=0$ on $[2s,\infty)$ and $-\frac{c}{t}\leq\phi'(t)\leq0$ on $(s,2s)$ for some constant $c>0$.
Define $D(r):=\{x\in M|f(x)\leq r\}$.

By Lemma 4.2, we have
\begin{eqnarray}
-\int_MR_{ijkl}\nabla_i\nabla_kR_{jl}\phi^2(f)e^{-f}&=&\int_M(|\nabla Ric|^2-\nabla_lR_{kj}\nabla_kR_{jl})(\phi^2)'e^{-f}\notag\\
&=&2\int_M(|\nabla Ric|^2-\nabla_lR_{kj}\nabla_kR_{jl})\phi\phi'e^{-f}\notag\\
&\leq&0.
\end{eqnarray}

Integrating (2.15) and using Lemma 4.1 and (4.27), we have
\begin{eqnarray}
&&\int_Mdiv^4Rm\phi^2(f)e^{-f}\notag\\
&=&\int_M\nabla_lR_{jk}\nabla_kR_{jl}\phi^2(f)e^{-f}-\int_M|\nabla Ric|^2\phi^2(f)e^{-f}-\int_MR_{ijkl}\nabla_i\nabla_kR_{jl}\phi^2(f)e^{-f}\notag\\
&\leq&c\int_M|Ric|^2|\nabla f|^2(\phi')^2e^{-f}+\frac{3}{4}\int_M|\nabla Ric|^2\phi^2(f)e^{-f}-\int_M|\nabla Ric|^2\phi^2(f)e^{-f}\notag\\
&\leq&\frac{c}{s^2}\int_{D(2s)\backslash D(s)}|Ric|^2|\nabla f|^2e^{-f}-\frac{1}{4}\int_M|\nabla Ric|^2\phi^2(f)e^{-f}.
\end{eqnarray}

It follows from Lemma 2.1, (2.12), (2.16) and Lemma 2.5 that
\[\int_M|Ric|^2|\nabla f|^2e^{-f}\leq\int_M|Ric|^2e^{-\alpha f}<+\infty\]
for some $\alpha\in(0,1]$. Therefore, $$\frac{c}{s^2}\int_{D(2s)\backslash D(s)}|Ric|^2|\nabla f|^2e^{-f}\rightarrow0$$ as $s\rightarrow+\infty$.

By taking $r\rightarrow+\infty$ in (4.28), we obtain $\int_M|\nabla Ric|^2e^{-f}=0$. Since $\int_M|\nabla Ric|^2e^{-f}<+\infty$, it follows from (2.20) that
\[\int_M|divRm|^2e^{-f}=\int_M|\nabla Ric|^2e^{-f}=0.\]

Hence,  $|divRm|=|\nabla Ric|=0$ $a. e$. .  Note that any gradient shrinking Ricci soliton is analytic in harmonic coordinates, we have $|divRm|=|\nabla Ric|=0$ on $M$.

It is clear $divRm=0$ implies that $M^n$ is radially flat.

By direct computation, we have
\[0\leq|\nabla Ric-\frac{\nabla R}{n}g|^2=|\nabla Ric|^2-\frac{|\nabla R|^2}{n}=-\frac{|\nabla R|^2}{n}\leq0.\]

Therefore, $R$ is a constant on $M$.

Since $M^n$ is radially flat and has constant scalar curvature, it follows from Lemma 2.4 that $(M^n,g)$ is rigid.$\hfill\Box$

Theorem 1.1 follows by combining Theorem 3.1 and Theorem 4.1.
\section{The proof of Theorem 1.2}
\renewcommand{\theequation}{\thesection.\arabic{equation}}
\indent\par
In this section, we give a direct proof of Theorem 1.2.
\\\\\textbf{Theorem 5.1} Let $(M^n,g)$ be a gradient shrinking Ricci soliton with (1.1). If $div^3Rm(\nabla f)=0$, then $(M^n,g)$ is rigid.
\\\\\textbf{Proof.} From (2.14), we have
\begin{eqnarray*}
div^3Rm(\nabla f)&=&\nabla_k\nabla_j\nabla_lR_{ijkl}\nabla_if\notag\\
&=&-R_{ijkl}\nabla_kR_{jl}\nabla_if\notag\\
&=&\frac{1}{2}(\nabla_iR_{ijkl})(\nabla_lR_{jk}-\nabla_kR_{jl})\notag\\
&=&-\frac{1}{2}|divRm|^2,
\end{eqnarray*}
where we used (2.4) in the third equality and (2.2) in the last.

Since $div^3Rm(\nabla f)=0$, $divRm=0$. It follows that $M$ is radially flat. Moreover, we have
\[\nabla_iR=2\nabla_lR_{il}=-2g^{jk}\nabla_lR_{ijkl}=0,\]
i.e. $R$ is a constant on $M$.

Since $M^n$ is radially flat and has constant scalar curvature, it follows from Lemma 2.4 that $(M^n,g)$ is rigid.$\hfill\Box$
\section{Under the condition of Weyl tensor}
\renewcommand{\theequation}{\thesection.\arabic{equation}}
\indent\par
In this section, we prove Theorems 1.3. The first step is to obtain the following formulas.
\\\\\textbf{Proposition 6.1} Let $(M^n,g)$ $(n\geq3)$ be a gradient shrinking Ricci soliton with (1.1), we have the following identities.
\begin{equation}
(divW)_{ijk}=\frac{n-3}{n-2}(divRm)_{ijk}-\frac{n-3}{2(n-1)(n-2)}(g_{ik}\nabla_jR-g_{jk}\nabla_iR),
\end{equation}
\begin{equation}
(div^2W)_{ik}=\frac{n-3}{n-2}(div^2Rm)_{ik}-\frac{n-3}{2(n-1)(n-2)}(g_{ik}\Delta R-\nabla_k\nabla_iR),
\end{equation}
\begin{equation}
(div^3W)_{i}=\frac{n-3}{n-2}(div^3Rm)_{i}+\frac{n-3}{2(n-1)(n-2)}R_{ik}\nabla_kR,
\end{equation}
and
\begin{equation}
div^4W=\frac{n-3}{n-2}div^4Rm+\frac{n-3}{2(n-1)(n-2)}(\frac{1}{2}|\nabla R|^2+R_{ik}\nabla_i\nabla_kR).
\end{equation}
\\\\\textbf{Proof.} By direct computation,
\begin{eqnarray*}
(divW)_{ijk}&=&\nabla_lW_{ijkl}\notag\\
&=&\nabla_lR_{ijkl}-\frac{1}{n-2}(g_{ik}\nabla_lR_{jl}-\nabla_iR_{jk}-g_{jk}\nabla_lR_{il}+\nabla_jR_{ik})\notag\\
&&+\frac{1}{(n-1)(n-2)}(g_{ik}\nabla_jR-g_{jk}\nabla_iR)\notag\\
&=&\nabla_lR_{ijkl}-\frac{1}{n-2}\nabla_lR_{ijkl}\notag\\
&&-\frac{1}{2(n-2)}(g_{ik}\nabla_jR-g_{jk}\nabla_iR)+\frac{1}{(n-1)(n-2)}(g_{ik}\nabla_jR-g_{jk}\nabla_iR)\notag\\
&=&\frac{n-3}{n-2}\nabla_lR_{ijkl}-\frac{n-3}{2(n-1)(n-2)}(g_{ik}\nabla_jR-g_{jk}\nabla_iR),
\end{eqnarray*}
where we used (2.8)in the second equality.

It follows from (6.29) that
\begin{eqnarray*}
&&(div^2W)_{ik}\notag\\
&=&\nabla_j\nabla_lW_{ijkl}\notag\\
&=&\frac{n-3}{n-2}\nabla_j\nabla_lR_{ijkl}-\frac{n-3}{2(n-1)(n-2)}(g_{ik}\Delta R-\nabla_k\nabla_iR),
\end{eqnarray*}

By (6.30), we have
\begin{eqnarray*}
(div^3W)_i&=&\nabla_k\nabla_j\nabla_lW_{ijkl}\notag\\
&=&\frac{n-3}{n-2}\nabla_k\nabla_j\nabla_lR_{ijkl}-\frac{n-3}{2(n-1)(n-2)}(\nabla_i\Delta R-\nabla_k\nabla_k\nabla_iR)\notag\\
&=&\frac{n-3}{n-2}\nabla_k\nabla_j\nabla_lR_{ijkl}+\frac{n-3}{2(n-1)(n-2)}R_{ik}\nabla_kR,
\end{eqnarray*}

From (6.31), we have
\begin{eqnarray*}
div^4W&=&\nabla_i\nabla_k\nabla_j\nabla_lW_{ijkl}\notag\\
&=&\frac{n-3}{n-2}\nabla_i\nabla_k\nabla_j\nabla_lR_{ijkl}+\frac{n-3}{2(n-1)(n-2)}\nabla_i(R_{ik}\nabla_kR)\notag\\
&=&\frac{n-3}{n-2}\nabla_i\nabla_k\nabla_j\nabla_lR_{ijkl}\notag\\
&&+\frac{n-3}{2(n-1)(n-2)}(\frac{|\nabla R|^2}{2}+R_{ik}\nabla_i\nabla_kR),
\end{eqnarray*}
$\hfill\Box$

As a corollary of Proposition 6.1, we have
\\\\\textbf{Corollary 6.1} Let $(M^n,g)$ $(n\geq3)$ be a gradient shrinking Ricci soliton with (1.1), we have the following identities.
\begin{equation}
(divW)_{ijk}=\frac{n-3}{n-2}(\nabla_jR_{ik}-\nabla_iR_{jk})-\frac{n-3}{2(n-1)(n-2)}(g_{ik}\nabla_jR-g_{jk}\nabla_iR),
\end{equation}
\begin{eqnarray}
(div^2W)_{ik}&=&\frac{n-3}{n-2}(2\lambda R_{ik}+\nabla_{\nabla f}R_{ik}-R_{ik}^2-R_{ijkl}R_{jl})-\frac{n-3}{2(n-1)}\nabla_i\nabla_kR\notag\\
&&-\frac{n-3}{2(n-1)(n-2)}(\nabla_{\nabla f}R+2\lambda R-2|Ric|^2)g_{ik},
\end{eqnarray}
\begin{equation}
(div^3W)_{i}=-\frac{n-3}{n-2}R_{ijkl}\nabla_kR_{jl}+\frac{n-3}{2(n-1)(n-2)}R_{ik}\nabla_kR,
\end{equation}
and
\begin{eqnarray}
div^4W&=&\frac{n-3}{n-2}(\nabla_lR_{jk}\nabla_kR_{jl}-|\nabla Ric|^2-R_{ijkl}\nabla_i\nabla_kR_{jl})\notag\\
&&+\frac{n-3}{2(n-1)(n-2)}(\frac{1}{2}|\nabla R|^2+R_{ik}\nabla_i\nabla_kR).
\end{eqnarray}
\\\\\textbf{Proof.} Applying (2.2) to (6.29), we obtain (6.33).

Applying (2.10) and (2.13) to (6.30), we can get (6.34).

Applying (2.14) to (6.31), we have (6.35).

Applying (2.15) to (6.32), we have (6.36).$\hfill\Box$
\\\\Next, we prove that a gradient shrinking Ricci soliton with $div^3W(\nabla f)=0$ is rigid.
\\\\\textbf{Theorem 6.1} Let $(M^n,g)$ be a gradient shrinking Ricci soliton with (1.1). If $div^3W(\nabla f)=0$, then $(M^n,g)$ is rigid.
\\\\\textbf{Proof.} By (6.35), we have
\begin{eqnarray}
&&div^3W(\nabla f)\notag\\
&=&-\frac{n-3}{n-2}R_{ijkl}\nabla_kR_{jl}\nabla_if+\frac{n-3}{2(n-1)(n-2)}R_{ik}\nabla_kR\nabla_if\notag\\
&=&\frac{n-3}{2(n-2)}(\nabla_iR_{ijkl})(\nabla_lR_{jk}-\nabla_kR_{jl})+\frac{n-3}{4(n-1)(n-2)}|\nabla R|^2\notag\\
&=&-\frac{n-3}{2(n-2)}|divRm|^2+\frac{n-3}{4(n-1)(n-2)}|\nabla R|^2,
\end{eqnarray}
where we used (2.4) and (2.8) in the second equality and (2.2) in the last.

It follows from (2.8) that $|\nabla R|^2\leq4|Ric|^2|\nabla f|^2$. By Lemma 2.1, (2.12), (2.16) and Lemma 2.5, we have
\begin{equation}
\int_M|\nabla R|^2e^{-f}\leq 4\int_M|Ric|^2|\nabla f|^2e^{-f}\leq\int_M|Ric|^2e^{-\alpha f}<+\infty,\notag
\end{equation}
where for some constant $\alpha\in(0,1]$.

Integrating (6.37) and using the condition of $div^3W(\nabla f)=0$, we obtain
\begin{equation}
\int_M|divRm|^2e^{-f}=\frac{1}{2(n-1)}\int_M|\nabla R|^2e^{-f}<+\infty.\notag
\end{equation}
It follows from (2.20) that
\begin{eqnarray}
\int_M|\nabla Ric|^2e^{-f}&=&\int_M|divRm|^2e^{-f}\notag\\
&=&\frac{1}{2(n-1)}\int_M|\nabla R|^2e^{-f}\notag\\
&\leq&\frac{n}{2(n-1)}\int_M|\nabla Ric|^2e^{-f},
\end{eqnarray}
where we used $|\nabla R|^2\leq n|\nabla Ric|^2$.

Note that $\frac{n}{2(n-1)}<1$, we conclude from (6.38) that
\begin{equation}
\int_M|divRm|^2e^{-f}=\int_M|\nabla R|^2e^{-f}=0,\notag
\end{equation}
i.e., $|divRm|=|\nabla R|=0$ $a. e.$ .

Since any gradient shrinking Ricci soliton is analytic in harmonic coordinates, $|divRm|=0$ on $M$. It follows that $M$ is radially flat. Moreover, $|\nabla R|=0$ on $M$, i.e., $R$ is a constant on $M$. By Lemma 2.4, $(M^n,g)$ is rigid.$\hfill\Box$
\section{Four-dimensional Case}
\renewcommand{\theequation}{\thesection.\arabic{equation}}
\indent\par
We prove Theorems 1.4 to 1.6 in this section. From Theorems 1.1 to 1.3, we only need to show the following classification theorem.
\\\\\textbf{Theorem 7.1} Let $(M^4,g)$ be a $4$-dimensional rigid gradient shrinking Ricci soliton with (1.1), then $(M^4,g)$ is either

(\rmnum{1}) Einstein, or

(\rmnum{2}) a finite quotient of the Gaussian shrinking soliton
$\mathbb{R}^4$, $\mathbb{R}^2\times\mathbb{S}^2$ or the round cylinder $\mathbb{R}\times\mathbb{S}^3$.

Before we prove Theorem 7.1, we present some results that are needed in the proof of Theorem 7.1.
\\\\\textbf{Lemma 7.1} (M. Fern\'{a}ndez-L\'{o}pez and E. Garc\'{i}a-R\'{i}o [9]) Let $(M^n,g,f)$ be an $n$-dimensional gradient shrinking Ricci soliton with constant scalar curvature, then $R\in\{0,\lambda,\cdot,\cdot\cdot,(n-1)\lambda,n\lambda\}$.
\\\\\textbf{Lemma 7.2} (M. Fern\'{a}ndez-L\'{o}pez and E. Garc\'{i}a-R\'{i}o [9]) No complete gradient shrinking Ricci soliton may
exist with $R=\lambda$.

Now we are ready to prove Theorem 7.1.
\\\\\textbf{Proof of Theorem 7.1:}

Note that $(M^4,g)$ is rigid, i.e., it is a finite quotient of $\mathbb{R}^k\times N^{4-k}$, where $N$ is an Einstein manifold and $k\in\{0,1,2,3,4\}$. It follows that $M^4$ has constant scalar curvature. Moreover, Lemma 7.1 and Lemma 7.2 imply that $R\in\{0,2\lambda,3\lambda,4\lambda\}$.

We denote by $\{e_i\}_{i=1}^4$ a local orthonormal frame of $M^4$ with $e_1=\frac{\nabla f}{|\nabla f|}$. Moreover, We use $\{\alpha_i\}_{i=1}^4$ to represent eigenvalues of the Ricci tensor with corresponding orthonormal eigenvectors $\{e_i\}_{i=1}^4$.

In the following, we divide the arguments into four cases:

$\bullet$ Case 1: $R\equiv0$. In this case, $(M^4,g,f)$ is a finite quotient of the Gaussian soliton $\mathbb{R}^4$.

$\bullet$ Case 2: $R\equiv2\lambda$. In this case, we have
\[(\alpha_1,\alpha_2,\alpha_3,\alpha_4)\in\{(\frac{\lambda}{2},\frac{\lambda}{2},\frac{\lambda}{2},\frac{\lambda}{2}),(0,\frac{2\lambda}{3},\frac{2\lambda}{3},\frac{2\lambda}{3}),(0,0,\lambda,\lambda),(0,0,0,2\lambda)\}.\]

It follows from (2.10) that $|Ric|^2=\lambda R=2\lambda^2$. Therefore, $(\alpha_1,\alpha_2,\alpha_3,\alpha_4)=(0,0,\lambda,\lambda)$. The rigidity of $(M^4,g)$ implies that it is a finite quotient of $\mathbb{R}^2\times N^2$ with positive scalar curvature. It is clear that $N^2$ has to be $\mathbb{S}^2$. Therefore, $(M^4,g)$ is a finite quotient of $\mathbb{R}^2\times\mathbb{S}^2$.

$\bullet$ Case 3: $R\equiv3\lambda$. In this case, we have
\[(\alpha_1,\alpha_2,\alpha_3,\alpha_4)\in\{(\frac{3\lambda}{4},\frac{3\lambda}{4},\frac{3\lambda}{4},\frac{3\lambda}{4}),(0,\lambda,\lambda,\lambda),(0,0,\frac{3\lambda}{2},\frac{3\lambda}{2}),(0,0,0,3\lambda)\}.\]

It follows from (2.10) that $|Ric|^2=\lambda R=3\lambda^2$. Therefore, $(\alpha_1,\alpha_2,\alpha_3,\alpha_4)=(0,\lambda,\lambda,\lambda)$.  The rigidity of $(M^4,g)$ implies that it is a finite quotient of $\mathbb{R}\times N^3$, where $N^3$ is Einstein with positive scalar curvature. It is clear that $N^3$ has to be $\mathbb{S}^3$. Therefore, $(M^4,g)$ is a finite quotient of $\mathbb{R}\times\mathbb{S}^3$.

$\bullet$ Case 4: $R\equiv4\lambda$. In this case, $(\alpha_1,\alpha_2,\alpha_3,\alpha_4)=(\lambda,\lambda,\lambda,\lambda)$, i.e., $(M^4,g)$ is Einstein with $Ric=\lambda g$.

We conclude that $(M^4,g)$ is either Einstein or a finite quotient of $\mathbb{R}^4$, $\mathbb{R}^2\times\mathbb{S}^2$ or $\mathbb{R}\times\mathbb{S}^3$.$\hfill\Box$
\section{Appendix}
\renewcommand{\theequation}{\thesection.\arabic{equation}}
\indent\par
G. Catino, P. Mastrolia and D. D. Monticelli [4] defined the fourth order divergence of Weyl tensor $div^4W$ to be $\nabla_k\nabla_j\nabla_l\nabla_iW_{ikjl}$. Moreover, they proved that a gradient shrinking Ricci soliton with $div^4W=0$ is rigid. It is clear from their proof that this result holds for $\nabla_k\nabla_j\nabla_l\nabla_iW_{ikjl}\leq0$.
\\\\\textbf{Remark 8.1} The definition of $div^4W$ in G. Catino, P. Mastrolia and D. D. Monticelli [4] differs from ours by a minus sign. To be more precise, we have
\begin{equation}
\nabla_k\nabla_j\nabla_l\nabla_iW_{ikjl}=\nabla_j\nabla_k\nabla_l\nabla_iW_{ijkl}=-\nabla_j\nabla_k\nabla_l\nabla_iW_{jikl}=-\nabla_i\nabla_k\nabla_l\nabla_jW_{ijkl}.
\end{equation}
It follows from (6.30) that $\nabla_j\nabla_lW_{ijkl}$ is symmetric on $i$ and $k$, then it is also symmetric on $j$ and $l$, i.e.,
\begin{equation}
\nabla_j\nabla_lW_{ijkl}=\nabla_l\nabla_jW_{ijkl}.
\end{equation}
Combining (8.39) and (8.40), we have
\[\nabla_k\nabla_j\nabla_l\nabla_iW_{ikjl}=-\nabla_i\nabla_k\nabla_j\nabla_lW_{ijkl}.\]
It is clear from (6.32) that
\[div^4W=\frac{n-3}{n-2}div^4Rm+\frac{n-3}{2(n-1)(n-2)}(\frac{1}{2}|\nabla R|^2+R_{ik}\nabla_i\nabla_kR).\]

The following theorems were proved by G. Catino, P. Mastrolia and D. D. Monticelli [4], we give a different proof here.
\\\\\textbf{Theorem 8.1} Let $(M^n,g)$ be a compact gradient shrinking Ricci soliton with (1.1). If $div^4W=0$, then $(M^n,g)$ is Einstein.
\\\\\textbf{Proof.} Integrating (6.36), we have
\begin{eqnarray}
&&\int_Mdiv^4We^{-f}\notag\\
&=&\frac{n-3}{n-2}\int_M(\nabla_lR_{jk}\nabla_kR_{jl}-|\nabla Ric|^2-R_{ijkl}\nabla_i\nabla_kR_{jl})e^{-f}\notag\\
&&+\frac{n-3}{2(n-1)(n-2)}\int_M(\frac{1}{2}|\nabla R|^2+R_{ik}\nabla_i\nabla_kR)e^{-f}\notag\\
&=&-\frac{n-3}{2(n-2)}\int_M|\nabla Ric|^2e^{-f}+\frac{n-3}{4(n-1)(n-2)}\int_M|\nabla R|^2e^{-f}\notag\\
&\leq&-\frac{n-3}{4n(n-1)}\int_M|\nabla R|^2e^{-f},
\end{eqnarray}
where we used Lemma 3.1, (2.5) and (2.7) in the second equality. Moreover, we used $|\nabla R|^2\leq n|\nabla Ric|^2$ in the inequality.

Since $div^4W=0$, it follows from (8.41) that  $\nabla R=0$ $a. e.$ . Note that any gradient shrinking Ricci soliton is analytic in harmonic coordinates, we have $\nabla R=0$ on $M$, i.e., $R$ is a constant on $M$. Therefore, $Ric(\nabla f,\nabla f)=\frac{1}{2}\langle\nabla R,\nabla f\rangle=0$. By Lemma 2.3, $(M^n,g)$ is Einstein.$\hfill\Box$
\\\\\textbf{Theorem 8.2} Let $(M^n,g)$ be a complete non-compact gradient shrinking Ricci soliton with (1.1). If $div^4W=0$, then $(M^n,g)$ is rigid.
\\\\\textbf{Proof.} Let $\phi:\mathbb{R_+}\rightarrow\mathbb{R}$ be a $C^2$ function with $\phi=1$ on $(0,s]$, $\phi=0$ on $[2s,\infty)$ and $-\frac{c}{t}\leq\phi'(t)\leq0$ on $(s,2s)$ for some constant $c>0$.
Define $D(r):=\{x\in M|f(x)\leq r\}$.

Integrating (6.37) we have
\begin{eqnarray}
&&\int_Mdiv^4W\phi^2(f)e^{-f}\notag\\
&=&\frac{n-3}{n-2}\int_M(\nabla_lR_{jk}\nabla_kR_{jl}-|\nabla Ric|^2-R_{ijkl}\nabla_i\nabla_kR_{jl})\phi^2(f)e^{-f}\notag\\
&&+\frac{n-3}{2(n-1)(n-2)}\int_M(\frac{1}{2}|\nabla R|^2+R_{ik}\nabla_i\nabla_kR)\phi^2(f)e^{-f}\notag\\
&\leq&c\int_{D(2s)\backslash D(s)}|Ric|^2|\nabla f|^2(\phi')^2e^{-f}-\frac{n-3}{4(n-2)}\int_M|\nabla Ric|^2\phi^2(f)e^{-f}\notag\\
&&+\frac{n-3}{4(n-1)(n-2)}\int_M|\nabla R|^2\phi^2e^{-f}-\frac{n-3}{(n-1)(n-2)}\int_MR_{ik}\nabla_kR\nabla_if\phi\phi'e^{-f}\notag\\
&\leq&\frac{c}{s^2}\int_{D(2s)\backslash D(s)}|Ric|^2|\nabla f|^2e^{-f}+\frac{n-3}{3(n-1)(n-2)}\int_M|\nabla R|^2\phi^2(f)e^{-f}\notag\\
&&-\frac{n-3}{2(n-2)}\int_M|\nabla Ric|^2\phi^2(f)e^{-f},
\end{eqnarray}
where we used Lemma 4.1 and Lemma 4.2 in the first inequality.

Applying $div^4W=0$ and $|\nabla Ric|\geq\frac{|\nabla R|^2}{n}$ to (8.42), we obtain
\begin{eqnarray}
0&\leq&\frac{c}{s^2}\int_{D(2s)\backslash D(s)}|Ric|^2|\nabla f|^2e^{-f}-\frac{(n-3)^2}{6n(n-1)(n-2)}\int_M|\nabla R|^2\phi^2(f)e^{-f}\notag\\
\end{eqnarray}

It follows from Lemma 2.1, (2.12), (2.16) and Lemma 2.5 that
\[\int_M|Ric|^2|\nabla f|^2e^{-f}\leq\int_M|Ric|^2e^{-\alpha f}<+\infty\]
for some $\alpha\in(0,1]$. Therefore, $$\frac{c}{s^2}\int_{D(2s)\backslash D(s)}|Ric|^2|\nabla f|^2e^{-f}\rightarrow0$$ as $s\rightarrow+\infty$.

By taking $r\rightarrow+\infty$ in (8.43), we obtain $\int_M|\nabla R|^2e^{-f}=0$. It follows that  $\nabla R=0$ $a. e.$ . Note that any gradient shrinking Ricci soliton is analytic in harmonic coordinates, we have $\nabla R=0$ on $M$, i.e., $R$ is a constant on $M$.

By taking $r\rightarrow+\infty$ in (8.43) and using  $div^4W=0$ and $|\nabla R|=0$, we obtain $\int_M|\nabla Ric|^2e^{-f}=0$. Since $\int_M|\nabla Ric|^2e^{-f}<+\infty$, it follows (2.20) that
\begin{equation}
\int_M|divRm|^2e^{-f}=\int_M|\nabla Ric|^2e^{-f}=0.
\end{equation}

Hence,  $|divRm|=0$ $a. e$. .  Note that any gradient shrinking Ricci soliton is analytic in harmonic coordinates, we have $|divRm|=0$ on $M$.

It is clear $divRm=0$ implies that $M^n$ is radially flat.

Since $M^n$ is radially flat and has constant scalar curvature, it follows from Lemma 2.4 that $(M^n,g)$ is rigid.$\hfill\Box$

From Theorem 7.1, Theorem 8.1 and Theorem 8.2, we have a classification theorem of $4$-dimensional gradient shrinking Ricci solitons with $div^4W=0$:
\\\\\textbf{Theorem 8.3} Let $(M^4,g)$ be a $4$-dimensional gradient shrinking Ricci soliton with (1.1). If $div^4W=0$, then $(M^4,g)$ is either

(\rmnum{1}) Einstein, or

(\rmnum{2}) a finite quotient of the Gaussian shrinking soliton
$\mathbb{R}^4$, $\mathbb{R}^2\times\mathbb{S}^2$ or the round cylinder $\mathbb{R}\times\mathbb{S}^3$.
\\\\\textbf{Remark 8.2} It is clear from the proof that Theorems 8.1 to 8.3 hold for $div^4W\geq0$. Moreover, it follows from (8.39) that Theorem 8.1 to 8.3 still hold if indices of $div^4W$ permutate.
\section*{Acknowledgements}
\renewcommand{\theequation}{\thesection.\arabic{equation}}
\indent\par
We would like to thank Professor Huai-Dong Cao for his encouragement and suggestions in improving the paper. The first author also thanks Professor Huai-Dong Cao for kindly invitation and warm hospitality during his stay at Lehigh University.


\begin{thebibliography}{99}
\bibitem{ref1} H. D. Cao, Q. Chen. On locally conformally flat gradient steady Ricci solitons. Transactions of the American Mathematical Society, 2012, 364: 2377-2391.
\bibitem{ref2} H. D. Cao, Q. Chen. On Bach-flat gradient shrinking Ricci solitons. Duke Mathematical Journal, 2013, 162(6): 1149-1169.
\bibitem{ref3} H. D. Cao, D. Zhou. On complete gradient shrinking Ricci solitons. Journal of Differential Geometry, 2010, 85(2): 175-186.
\bibitem{red4} G. Catino, P. Mastrolia, D. D. Monticelli, Gradient Ricci solitons with vanishing conditions on Weyl, arXiv: 1602.00534v2 [math.DG].
\bibitem{ref5} B. L. Chen. Strong uniqueness of the Ricci flow. Journal of Differential Geometry, 2009, 86(2): 362-382.
\bibitem{ref6} X. Chen, Y. Wang. On four-dimensional anti-self-dual gradient Ricci solitons. Journal of Geometric Analysis, 2015, 25: 1335-1343.
\bibitem{ref7} M. Eminenti, G. La Nave, C. Mantegazza. Ricci solitons: the equation point of view. Manuscripta Mathematica, 2008, 127: 345-367.
\bibitem{ref8} M. Fern\'{a}ndez-L\'{o}pez, E. Garc\'{i}a-R\'{i}o. Rigidity of shrinking Ricci solitons. Mathematische Zeitschrift, 2011, 269(1): 461-466.
\bibitem{ref9} M. Fern\'{a}ndez-L\'{o}pez, E. Garc\'{i}a-R\'{i}o. On gradient Ricci solitons with constant scalar curvature. Proceedings of the American Mathematical Society, 2016, 144: 369-378.
\bibitem{ref10} R. S. Hamilton. The formation of singularities in the Ricci flow, Surveys in Differential
Geometry (Cambridge, MA, 1993), 2: 7-136.
\bibitem{ref11} O. Munteanu, N. Sesum. On gradient Ricci solitons. Journal of Geometric Analysis, 2013, 23: 539-561.
\bibitem{ref12} A. Naber. Noncompact shrinking 4-solitons with nonnegative curvature. Journal f\"{u}r die reine und angewandte Mathematik, 2007, 645(2): 125-153.
\bibitem{ref13} P. Petersen, W. Wylie. Rigidity of gradient Ricci solitons. Pacific Journal of Mathematics, 2009, 241(2): 329-345.
\bibitem{ref14} P. Petersen, W. Wylie. On the classification of gradient Ricci solitons. Geometry \& Topology, 2010, 14(4): 2277-2300.
\bibitem{ref15} S. Tachibana. A theorem of Riemannian manifolds of positive curvature operator. Proceedings of the Japan Academy, 1974, 50(4): 301-302.
\bibitem{ref16} V. Timofte. On the positivity of symmetric polynomial functions. Part \Rmnum{1}: general results. Journal of Mathematical Analysis and Applications, 2003, 284(1):174-190.
\bibitem{ref17} J. Y. Wu, P. Wu, W. Wylie. Gradient shrinking Ricci solitons of half harmonic Weyl curvature. arXiv: 1410.7303v1 [math.DG].
\bibitem{ref18} Z. H. Zhang. Gadient shrinking solitons with vanishing Weyl tensor. Pacific Journal of Mathematics, 2009, 242(1): 189-200.
\end{thebibliography}
\end{document}